\documentclass[10.5pt]{article}
\usepackage{amsmath}
\usepackage{amssymb}
\usepackage{amscd}
\usepackage{xspace}
\usepackage{verbatim}
\newcommand{\mysection}[1]{
\section{#1}\setcounter{equation}{0}}
\title{\bf On uniqueness of large solutions of nonlinear parabolic equations in nonsmooth domains}
\author{{\bf Waad Al Sayed}\quad
 {\bf Laurent V\'eron}\\[2mm]
{\small Laboratoire de Math\'ematiques et Physique Th\'eorique, }\\
{\small  Universit\'e Fran\c{c}ois Rabelais,  Tours,  FRANCE}}

\date{}
\begin{document}
\maketitle
{\small {\bf Abstract} We study the existence and uniqueness of the positive solutions of the problem (P):  $\partial_tu-\Delta u+u^q=0$ ($q>1$) in $\Omega\times (0,\infty)$, $u=\infty$ on $\partial\Omega\times (0,\infty)$ and $u(.,0)\in L^1(\Omega)$, when $\Omega$ is a bounded domain in $\mathbb R^N$.  We construct a maximal solution, prove that this maximal solution is a large solution whenever $q<N/(N-2)$ and it is unique if $\partial\Omega=\partial\overline\Omega^c$. If $\partial\Omega$ has the local graph property, we prove that there exists at most one solution to problem (P).
}

\vspace{1mm}
\hspace{.05in}

\newcommand{\txt}[1]{\;\text{ #1 }\;}
\newcommand{\tbf}{\textbf}
\newcommand{\tit}{\textit}
\newcommand{\tsc}{\textsc}
\newcommand{\trm}{\textrm}
\newcommand{\mbf}{\mathbf}
\newcommand{\mrm}{\mathrm}
\newcommand{\bsym}{\boldsymbol}
\newcommand{\scs}{\scriptstyle}
\newcommand{\sss}{\scriptscriptstyle}
\newcommand{\txts}{\textstyle}
\newcommand{\dsps}{\displaystyle}
\newcommand{\fnz}{\footnotesize}
\newcommand{\scz}{\scriptsize}
\newcommand{\be}{
\begin{equation}
}
\newcommand{\bel}[1]{
\begin{equation}
\label{#1}}
\newcommand{\ee}{
\end{equation}
}
\newcommand{\eqnl}[2]{
\begin{equation}
\label{#1}{#2}
\end{equation}
}
\newtheorem{subn}{\name}
\renewcommand{\thesubn}{}
\newcommand{\bsn}[1]{\def\name{#1}
\begin{subn}}
\newcommand{\esn}{
\end{subn}}
\newtheorem{sub}{\name}[section]
\newcommand{\dn}[1]{\def\name{#1}}   
\newcommand{\bs}{
\begin{sub}}
\newcommand{\es}{
\end{sub}}
\newcommand{\bsl}[1]{
\begin{sub}\label{#1}}
\newcommand{\bth}[1]{\def\name{Theorem}
\begin{sub}\label{t:#1}}
\newcommand{\blemma}[1]{\def\name{Lemma}
\begin{sub}\label{l:#1}}
\newcommand{\bcor}[1]{\def\name{Corollary}
\begin{sub}\label{c:#1}}
\newcommand{\bdef}[1]{\def\name{Definition}
\begin{sub}\label{d:#1}}
\newcommand{\bprop}[1]{\def\name{Proposition}
\begin{sub}\label{p:#1}}
\newcommand{\R}{\eqref}
\newcommand{\rth}[1]{Theorem~\ref{t:#1}}
\newcommand{\rlemma}[1]{Lemma~\ref{l:#1}}
\newcommand{\rcor}[1]{Corollary~\ref{c:#1}}
\newcommand{\rdef}[1]{Definition~\ref{d:#1}}
\newcommand{\rprop}[1]{Proposition~\ref{p:#1}}
\newcommand{\BA}{
\begin{array}}
\newcommand{\EA}{
\end{array}}
\newcommand{\BAN}{\renewcommand{\arraystretch}{1.2}
\setlength{\arraycolsep}{2pt}
\begin{array}}
\newcommand{\BAV}[2]{\renewcommand{\arraystretch}{#1}
\setlength{\arraycolsep}{#2}
\begin{array}}
\newcommand{\BSA}{
\begin{subarray}}
\newcommand{\ESA}{
\end{subarray}}
\newcommand{\BAL}{
\begin{aligned}}
\newcommand{\EAL}{
\end{aligned}}
\newcommand{\BALG}{
\begin{alignat}}
\newcommand{\EALG}{
\end{alignat}}
\newcommand{\BALGN}{
\begin{alignat*}}
\newcommand{\EALGN}{
\end{alignat*}}
\newcommand{\note}[1]{\textit{#1.}\hspace{2mm}}
\newcommand{\Proof}{\note{Proof}}
\newcommand{\qeda}{\hspace{10mm}\hfill $\square$}
\newcommand{\qed}{\\
${}$ \hfill $\square$}
\newcommand{\Remark}{\note{Remark}}
\newcommand{\modin}{$\,$\\
[-4mm] \indent}
\newcommand{\forevery}{\quad \forall}
\newcommand{\set}[1]{\{#1\}}
\newcommand{\setdef}[2]{\{\,#1:\,#2\,\}}
\newcommand{\setm}[2]{\{\,#1\mid #2\,\}}
\newcommand{\lra}{\longrightarrow}
\newcommand{\lla}{\longleftarrow}
\newcommand{\llra}{\longleftrightarrow}
\newcommand{\Lra}{\Longrightarrow}
\newcommand{\Lla}{\Longleftarrow}
\newcommand{\Llra}{\Longleftrightarrow}
\newcommand{\warrow}{\rightharpoonup}
\newcommand{
\paran}[1]{\left (#1 \right )}
\newcommand{\sqbr}[1]{\left [#1 \right ]}
\newcommand{\curlybr}[1]{\left \{#1 \right \}}
\newcommand{\abs}[1]{\left |#1\right |}
\newcommand{\norm}[1]{\left \|#1\right \|}
\newcommand{
\paranb}[1]{\big (#1 \big )}
\newcommand{\lsqbrb}[1]{\big [#1 \big ]}
\newcommand{\lcurlybrb}[1]{\big \{#1 \big \}}
\newcommand{\absb}[1]{\big |#1\big |}
\newcommand{\normb}[1]{\big \|#1\big \|}
\newcommand{
\paranB}[1]{\Big (#1 \Big )}
\newcommand{\absB}[1]{\Big |#1\Big |}
\newcommand{\normB}[1]{\Big \|#1\Big \|}

\newcommand{\thkl}{\rule[-.5mm]{.3mm}{3mm}}
\newcommand{\thknorm}[1]{\thkl #1 \thkl\,}
\newcommand{\trinorm}[1]{|\!|\!| #1 |\!|\!|\,}
\newcommand{\bang}[1]{\langle #1 \rangle}
\def\angb<#1>{\langle #1 \rangle}
\newcommand{\vstrut}[1]{\rule{0mm}{#1}}
\newcommand{\rec}[1]{\frac{1}{#1}}
\newcommand{\opname}[1]{\mbox{\rm #1}\,}
\newcommand{\supp}{\opname{supp}}
\newcommand{\dist}{\opname{dist}}
\newcommand{\myfrac}[2]{{\displaystyle \frac{#1}{#2} }}
\newcommand{\myint}[2]{{\displaystyle \int_{#1}^{#2}}}
\newcommand{\mysum}[2]{{\displaystyle \sum_{#1}^{#2}}}
\newcommand {\dint}{{\displaystyle \int\!\!\int}}
\newcommand{\q}{\quad}
\newcommand{\qq}{\qquad}
\newcommand{\hsp}[1]{\hspace{#1mm}}
\newcommand{\vsp}[1]{\vspace{#1mm}}
\newcommand{\ity}{\infty}
\newcommand{\prt}{\partial}
\newcommand{\sms}{\setminus}
\newcommand{\ems}{\emptyset}
\newcommand{\ti}{\times}
\newcommand{\pr}{^\prime}
\newcommand{\ppr}{^{\prime\prime}}
\newcommand{\tl}{\tilde}
\newcommand{\sbs}{\subset}
\newcommand{\sbeq}{\subseteq}
\newcommand{\nind}{\noindent}
\newcommand{\ind}{\indent}
\newcommand{\ovl}{\overline}
\newcommand{\unl}{\underline}
\newcommand{\nin}{\not\in}
\newcommand{\pfrac}[2]{\genfrac{(}{)}{}{}{#1}{#2}}

\def\ga{\alpha}     \def\gb{\beta}       \def\gg{\gamma}
\def\gc{\chi}       \def\gd{\delta}      \def\ge{\epsilon}
\def\gth{\theta}                         \def\vge{\varepsilon}
\def\gf{\phi}       \def\vgf{\varphi}    \def\gh{\eta}
\def\gi{\iota}      \def\gk{\kappa}      \def\gl{\lambda}
\def\gm{\mu}        \def\gn{\nu}         \def\gp{\pi}
\def\vgp{\varpi}    \def\gr{\rho}        \def\vgr{\varrho}
\def\gs{\sigma}     \def\vgs{\varsigma}  \def\gt{\tau}
\def\gu{\upsilon}   \def\gv{\vartheta}   \def\gw{\omega}
\def\gx{\xi}        \def\gy{\psi}        \def\gz{\zeta}
\def\Gg{\Gamma}     \def\Gd{\Delta}      \def\Gf{\Phi}
\def\Gth{\Theta}
\def\Gl{\Lambda}    \def\Gs{\Sigma}      \def\Gp{\Pi}
\def\Gw{\Omega}     \def\Gx{\Xi}         \def\Gy{\Psi}

\def\CS{{\mathcal S}}   \def\CM{{\mathcal M}}   \def\CN{{\mathcal N}}
\def\CR{{\mathcal R}}   \def\CO{{\mathcal O}}   \def\CP{{\mathcal P}}
\def\CA{{\mathcal A}}   \def\CB{{\mathcal B}}   \def\CC{{\mathcal C}}
\def\CD{{\mathcal D}}   \def\CE{{\mathcal E}}   \def\CF{{\mathcal F}}
\def\CG{{\mathcal G}}   \def\CH{{\mathcal H}}   \def\CI{{\mathcal I}}
\def\CJ{{\mathcal J}}   \def\CK{{\mathcal K}}   \def\CL{{\mathcal L}}
\def\CT{{\mathcal T}}   \def\CU{{\mathcal U}}   \def\CV{{\mathcal V}}
\def\CZ{{\mathcal Z}}   \def\CX{{\mathcal X}}   \def\CY{{\mathcal Y}}
\def\CW{{\mathcal W}} \def\CQ{{\mathcal Q}}
\def\BBA {\mathbb A}   \def\BBb {\mathbb B}    \def\BBC {\mathbb C}
\def\BBD {\mathbb D}   \def\BBE {\mathbb E}    \def\BBF {\mathbb F}
\def\BBG {\mathbb G}   \def\BBH {\mathbb H}    \def\BBI {\mathbb I}
\def\BBJ {\mathbb J}   \def\BBK {\mathbb K}    \def\BBL {\mathbb L}
\def\BBM {\mathbb M}   \def\BBN {\mathbb N}    \def\BBO {\mathbb O}
\def\BBP {\mathbb P}   \def\BBR {\mathbb R}    \def\BBS {\mathbb S}
\def\BBT {\mathbb T}   \def\BBU {\mathbb U}    \def\BBV {\mathbb V}
\def\BBW {\mathbb W}   \def\BBX {\mathbb X}    \def\BBY {\mathbb Y}
\def\BBZ {\mathbb Z}

\def\GTA {\mathfrak A}   \def\GTB {\mathfrak B}    \def\GTC {\mathfrak C}
\def\GTD {\mathfrak D}   \def\GTE {\mathfrak E}    \def\GTF {\mathfrak F}
\def\GTG {\mathfrak G}   \def\GTH {\mathfrak H}    \def\GTI {\mathfrak I}
\def\GTJ {\mathfrak J}   \def\GTK {\mathfrak K}    \def\GTL {\mathfrak L}
\def\GTM {\mathfrak M}   \def\GTN {\mathfrak N}    \def\GTO {\mathfrak O}
\def\GTP {\mathfrak P}   \def\GTR {\mathfrak R}    \def\GTS {\mathfrak S}
\def\GTT {\mathfrak T}   \def\GTU {\mathfrak U}    \def\GTV {\mathfrak V}
\def\GTW {\mathfrak W}   \def\GTX {\mathfrak X}    \def\GTY {\mathfrak Y}
\def\GTZ {\mathfrak Z}   \def\GTQ {\mathfrak Q}

\font\Sym= msam10 
\def\SYM#1{\hbox{\Sym #1}}
\newcommand{\bdw}{\prt\Gw\xspace}
\medskip
\mysection {Introduction}
Let $q>1$ and let  $\Gw$ be a bounded domain in $\BBR^N$ with boundary $\prt\Gw:=\Gg$. It has been proved by Keller \cite {Ke} and Osserman \cite {Oss} that there exists a {\it maximal solution} $\overline u$ to the stationnary equation
\begin{equation}\label{ell1}
-\Gd u+|u|^{q-1}u=0\quad\text{in }\Gw.
\end{equation} 
When $1<q<N/(N-2)$ this maximal solution is a {\it large solution} in the sense that
\begin{equation}\label{ell2}
\lim_{\gr(x)\to 0}\overline u(x)=\infty
\end{equation} 
where $\gr(x)=\dist (x,\prt\Gw)$. Furthermore V\'eron proves in \cite {Ve1} that $\overline u$ is the unique large solution whenever $\prt\Gw=\prt\overline\Gw^c$. When $q\geq N/(N-2)$  his proof of uniqueness does not apply. Marcus and V\'eron prove in \cite{MV1} that, there exists at most one large solution, provided $\prt\Gw$ is locally the graph of a continuous function. The aim of this article is to extend these questions to the parabolic equation
\begin{equation}\label{E1}
\prt_{t}u-\Gd u+|u|^{q-1}u=0\quad\text{in }\Gw\ti (0,\infty).
\end{equation} 
We are interested into positive solutions which satisfy 
\begin{equation}\label{T1}
\lim_{t\to 0}u(.,t)=f \quad\text{in }L^1_{loc}(\Gw),
\end{equation} 
where $f\in L^{1}_{loc\,+}(\Gw)$ and
\begin{equation}\label{T2}
\lim_{(x,t)\to (y,s)}u(x,t)=\infty \quad\forall (y,s)\in \Gg\ti (0,\infty).
\end{equation} 
Notice that if the initial and boundary conditions are exchanged, i.e. $u (.,t)$ blows-up when $t\to 0$ and coincides with a locally integrable function on  $\Gg\ti (0,\infty)$, this problem is associated with the study of the initial trace, and much work has been done by Marcus and V\'eron \cite {MV3} in the case of a smooth domain. In particular they obtain the existence and uniqueness when $q$ is subcritical, i.e. $1<q<1+2/N$. \medskip

In this article we prove two series of results:\medskip

\noindent{\bf Theorem A }{\it Assume $q>1$ and $\Gw$ is a bounded domain. Then for any $f\in L^{1}_{loc\,+}(\Gw)$ there exists a maximal solution $\overline u_{f}$ to problem (\ref{E1}) satisfying (\ref{T1}). If $1<q<N/(N-2)$, $\overline u_{f}$ satisfies (\ref{T2}). At end, if $1<q<N/(N-2)$ and $\prt\Gw=\prt\overline{\Gw}^c$, $\overline u_{f}$ is the unique solution of the problem which satisfies  (\ref{T2}).}\smallskip

The proof of uniqueness is based upon the construction of self-similar solutions of (\ref{E1}) in $\BBR^N\setminus\{0\}\ti (0,\infty)$, with a persistent strong singularity on the axis $\{0\}\ti (0,\infty)$ and a zero initial trace on $\BBR^N\setminus\{0\}$. This solution, which is studied in Appendix, is reminiscent of the very singular solution of Brezis, Peletier and Terman \cite {BPT}, although the method of construction is far different. The uniqueness is a delicate adaptation to the parabolic framework of the proof by contradiction of \cite{Ve1}.\medskip

\noindent{\bf Theorem B }{\it Assume $q>1$, $\Gw$ is a bounded domain and $\prt\Gw$, is locally a continuous graph. Then for any $f\in L^{1}_{loc\,+}(\Gw)$ there exists at most one solution to problem (\ref{E1}) satisfying (\ref{T1}) and (\ref{T2}).}\smallskip

For proving this result, we adapt the idea which was introduced in \cite{MV1} of constructing local super and subsolutions by small translations of the domain, but the non-uniformity of the boundary blow-up creates an extra-difficulty. In an appendix we study a self-similar equation which plays a key-role in our construction,
\begin{equation}\label{S1}\left\{\BA{l}
H''+\left(\myfrac{N-1}{r}+\myfrac{r}{2}\right)H'+\myfrac{1}{q-1}H-|H|^{q-1}=0\\[2mm]\phantom{\left(----\right)'+''\frac{1}{q-1}H'-|H|^{q-1}}\lim_{r\to 0}H(r)=\infty\\[2mm]\phantom{\left(\right)'+''\frac{1}{q-1}H'-|H|^{q-1}}
\lim_{r\to\infty}r^{2/(q-1)}H(r)=0.
\EA\right.\end{equation} 
We prove the existence and the uniqueness of the positive solution of (\ref{S1}) when $1<q<N/(N-2)$ and we give precise asymptotics when $r\to 0$ and $r\to \infty$.
\medskip

This article is organised as follows: 1- Introduction. 2- The maximal solution 3- The case $1<q<N/(N-2)$. 4- The local continuous graph property. 5- Appendix.
\section{The maximal solution}
\setcounter{equation}{0}

In this section $\Gw$ is an open domain of $\BBR^N$, with a compact boundary $\Gg:=\prt\Gw$. If $G$ is any open subset of $\BBR^N$ and $0<T\leq\infty$, we denote $Q_{T}^G:=G\ti (0,T)$. If $f\in L^1_{loc\,+}(\Gw)$, we consider the problem
\begin{equation}\label{M1}\left\{\BA {l}
\prt_{t}u-\Gd u+|u|^{q-1}u=0\quad\text{in }Q_{\infty}^\Gw\\[2mm]
\phantom{..........}
\lim_{t\to 0}u(.,t)=f(.)\quad\text{in }L^1_{loc}(\Gw)\\[2mm]
\lim_{(x,t)\to (y,s)}u(x,t)=\infty \quad\forall (y,s)\in \Gg\ti (0,\infty).
\EA\right.\end{equation}
By the next result, we reduce the lateral blow-up condition by a locally uniform one in which we set $\gr(x)=\dist(x,\Gg)$.
\blemma {equiv}The following two conditions are equivalent
\begin{equation}\label{C3}
\lim_{(x,t)\to (y,s)}u(x,t)=\infty \quad\forall (y,s)\in \Gg\ti (0,\infty)
\end{equation}
and
\begin{equation}\label{C4}\lim_{\gr(x)\to 0}u(x,t)=\infty \quad\text{uniformly on }[\gt,T],
\end{equation}
for any $0<\gt<T<\infty$.
\es
\Proof It is clear that (\ref{C4}) is equivalent to the fact that (\ref{C3}) holds uniformly on $\Gg\ti [\gt,T]$. By contradiction, we assume that (\ref{C3}) does not hold uniformly  for some $T>\gt>0$. Then there exists $\gb>0$ such that for any $\gd>0$, there exist two couples $(y_{\gd},s_{\gd})\in \Gg\ti [\gt,T]$ and $(x_{\gd},t_{\gd})\in \Gw\ti [\gt,T]$ such that 
\begin{equation}\label{C4'}|x_{\gd}-y_{\gd}|+|t_{\gd}-s_{\gd}|\leq\gd\quad\text{and }\;u(x_{\gd},t_{\gd})\leq\gb.
\end{equation}
Taking $\gd=1/n$, $n\in\BBN^{*}$, we can assume that $\{\gd\}$ is discrete and that 
$y_{\gd}\to y\in\Gg$ and $s_{\gd}\to s\in [\gt,T]$. Thus 
$x_{\gd}\to y$ and $t_{\gd}\to s$. Therefore
(\ref{C4'}) contradicts (\ref{C3}).\qeda
\bth{max}For any $q>1$ and $f\in L^1_{loc\,+}(\Gw)$, there exists a maximal solution $u:=\overline u_{f}$ of 
\begin{equation}\label{E1}
\prt_{t}u-\Gd u+|u|^{q-1}u=0\quad\text{in }Q_{\infty}^\Gw
\end{equation}
which satisfies 
\begin{equation}\label{C1}
\lim_{t\to 0}u(.,t)=f(.)\quad\text{in }L^1_{loc}(\Gw).
\end{equation}
\es
\Proof Let $\Gw_{n}$ be an increasing sequence of smooth bounded domains such that $\overline\Gw_{n}\subset\Gw_{n+1}\subsetÊ\Gw$ and $\cup\Gw_{n}=\Gw$. For each $n$ let $u_{n,f}$ be the increasing limit when $k\to\infty$ of the $u_{n,k,f}$ solution of
\begin{equation}\label {Ek}\left\{\BA {l}
\prt_{t}u_{n,k,f}-\Gd u_{n,k,f}+u_{n,k,f}^q=0\quad\text{in }Q^{\Gw_{n}}_{\infty}\\[2mm]
\phantom{,,,,,,,,}u_{n,k,f}(x,t)=k\quad\text{in }\prt\Gw_{n}\ti (0,\infty)\\[2mm]
\phantom{,,,,,,,,}u_{n,k,f}(x,0)=f\chi_{_{\Gw_{n}}}\quad\text{in }\Gw_{n}.
\EA\right.
\end{equation}
By the maximum principle and a standard approximation argument $n\mapsto u_{n,k,f}$ is decreasing thus $n\mapsto u_{n,f}$ too. The limit $\overline u_{f}$ of the  $u_{n,f}$ satisfies (\ref{E1}) and (\ref{C1}). It is independent of the exhaustion $\{\Gw_{n}\}$ of $\Gw$. Let $u$ be a positive  solution of (\ref{E1}) in $Q^\Gw_{\infty}$ which satisfies (\ref{C1}). Since the initial trace of $u$ is a locally integrable function, $u^q\in L^1_{loc}(\Gw\ti [0,\infty))$. By Fubini we can assume that, for any $n$, $u\in L^1_{loc}(\prt\Gw_{n}\ti [0,\infty))$. Because $(u-u_{n,k,f})_{+}\leq u$ and tends to $0$ when $k\to\infty$, it follows by Lebesgue's theorem that
$$\lim_{k\to\infty}\norm{(u-u_{n,k,f})_{+}}_{L^1(\prt\Gw_{n}\ti (0,T))}=0\quad\forall T>0.
$$
Applying the maximum principle in $\Gw_{n}\ti (0,\infty)$ yields to
$$u\leq \lim_{k\to\infty}u_{n,k,f}=u_{n,f}\Longrightarrow
u\leq \lim_{n\to\infty}u_{n,f}=\overline u_{f}.
$$
\qeda
\bth{min}For any $q>1$ and $f\in L^1_{loc\,+}(\Gw)$, there exists a minimal nonnegative solution $\underline u_{f}$ of (\ref{E1}) in $Q_{\infty}^\Gw$ which satisfies (\ref{C1}).
\es
\Proof The scheme of the construction is similar to the one of $\overline u_{f}$: with the same exhaustion $\{\Gw_{n}\}$ of $\Gw$, we consider
the solution $u_{n,0,f}$ solution of
\begin{equation}\label {E0}\left\{\BA {l}
\prt_{t}u_{n,0,f}-\Gd u_{n,0,f}+u_{n,0,f}^q=0\quad\text{in }Q^{\Gw_{n}}_{\infty}\\[2mm]
\phantom{,,,,,,,,}u_{n,0,f}(x,t)=0\quad\text{in }\prt\Gw_{n}\ti (0,\infty)\\[2mm]
\phantom{,,,,,,,,}u_{n,0,f}(x,0)=f\chi_{_{\Gw_{n}}}\quad\text{in }\Gw_{n}.
\EA\right.
\end{equation}
By the maximum principle, $n\mapsto u_{n,0,f}$ is increasing and dominated by $\overline u_{f}$. Therefore it converges to some solution 
$\underline u_{f}$ of (\ref{E1}), which satisfies (\ref{C1}) as $u_{n,0,f}$ and $\overline u_{f}$ do it. Using the same argument as in the proof of \rth{max}, there holds $u_{n,0,f}\leq u$ in $Q^{\Gw_{n}}_{\infty}$ for a suitable exhaustion. Thus $\underline u_{f}\leq u$.\qeda\medskip

\noindent\Remark Because of the lack of regularity of $\prt\Gw$, there is no reason for $\overline u_{f}$ (resp $\underline u_{f}$) to tend to infinity (resp. zero) on $\prt\Gw\ti (0,\infty)$. \medskip

The next statement will be very usefull for proving uniqueness results.

\bth {contraction}Assume $q>1$, $f\in L^1_{loc\,+}(\Gw)$ and $u_{f}$ is a nonnegative solution of (\ref{E1}) satisfying (\ref{C1}). Then there exists a nonnegative solution $u_{0}$ of (\ref{E1}) satisfying
\begin{equation}\label{D0}
\lim_{t\to 0}u_{0}(.,t)=0\quad\text{in }L^1_{loc}(\Gw),
\end{equation}
such that
\begin{equation}\label{D1}
0\leq u_{f}-\underline u_{f}\leq u_{0}\leq u_{f},
\end{equation}
and
\begin{equation}\label{D2}
0\leq\overline u_{f}-u_{f}\leq \overline u_{0}-u_{0}.
\end{equation}
\es
\Proof {\it Step 1: construction of $u_{0}$. }The function $w=u_{f}-\underline u_{f}$ is a nonnegative subsolution of (\ref{E1}) which satisfies
$$\lim_{t\to 0}w(.,t)=0\quad\text{in }L^1_{loc}(\Gw).
$$
Using the above considered exhaustion of $\Gw$, we denote by 
$v_{n}$ the solution of 
\begin{equation}\label {Em}\left\{\BA {l}
\prt_{t}v_{n}-\Gd v_{n}+v_{n}^q=0\quad\text{in }Q^{\Gw_{n}}_{\infty}\\[2mm]
\phantom{,,,,,,,,,}v_{n}(x,t)=u_{f}-\underline u_{f}\quad\text{in }\prt\Gw_{n}\ti (0,\infty)\\[2mm]
\phantom{,,,,,,,,,}v_{n}(x,0)=0\quad\text{in }\Gw_{n}.
\EA\right.
\end{equation}
By the maximum principle 
$$u_{f}-\underline u_{f}\leq v_{n}\leq u_{f}\quad\text{in }Q_{\infty}^{\Gw_{n}}.$$
Therefore $v_{n+1}\geq v_{n}$  on $\prt\Gw_{n}\ti (0,\infty)$; this implies that the same inequality holds in $Q_{\infty}^{\Gw_{n}}$. If we denote by $u_{0}$ the limit of the $\{v_{n}\}$, it is a solution of (\ref{E1}) in $Q_{\infty}^{\Gw}$. For any compact $K\in \Gw$, there exists $n_{K}$ and $\ga>0$ such that $\dist (K,\Gw_{n}^c)\geq \ga$ for $n\geq n_{K}$ therefore $v_{n}$ remains uniformly bounded on $K$ by Brezis-Friedman estimate \cite{BF}. Thus the local equicontinuity of the $v_{n}$ (consequence of the regularity theory for parabolic equations) implies that $u_{0}$ satisfies (\ref{D0}).\smallskip

\noindent{\it Step 2: proof of (\ref{D2}). }We follow a method introduced in \cite{MV2} in a different context. For $n\in\BBN$ and $k>0$ fixed, we  set 
$$Z_{f,n}=u_{f,n}-u_{f}\quad\text{and }\;Z_{0,n}= u_{0,n}-u_{0},
$$
where we assume that the $n$ are chosen such that $u_{f},u_{0}\in L^1_{loc}(\prt\Gw_{n}\ti [0,\infty))$, 
and
$$\gf(r,s)=\left\{\BA {l}\myfrac{r^q-s^q}{r-s}\quad\text{if }\;r\neq s\\[2mm]
0\qquad\quad\quad\text{if }\;r= s.
\EA\right.$$
By convexity, 
$$\left\{\BA {l}r_{0}\geq s_{0},\;r_{1}\geq s_{1}\\
r_{1}\geq r_{0},\;s_{1}\geq s_{0}
\EA\right.\Longrightarrow \gf(r_{1},s_{1})\geq \gf(r_{0},s_{0}).$$
Therefore
$$\gf(u_{f,n},u_{f})\geq \gf(u_{0,n},u_{0})\quad\text{in }Q^{\Gw_{n}}_{T},
$$
and
$$\BA {l}
0=\prt _{t}(Z_{f,n}-Z_{0,n})-\Gd (Z_{f,n}-Z_{0,n})
+u^q_{f,n}-u^q_{f}-u^q_{0,n}+u^q_{0}\\[2mm]\phantom{0}=
\prt _{t}(Z_{f,n}-Z_{0,n})-\Gd (Z_{f,n}-Z_{0,n})
+\gf(u_{f,n},u_{f})Z_{f,n}-\gf(u_{0,n},u_{0})Z_{0,n},
\EA$$
which implies
$$\prt _{t}(Z_{f,n}-Z_{0,n})-\Gd (Z_{f,n}-Z_{0,n})
+\gf(u_{f,n},u_{f})(Z_{f,n}-Z_{0,n})\leq 0.
$$
But $Z_{f,n}-Z_{0,n}=0$ in $\Gw_{n}\ti\{0\}$ and 
$$\myint{0}{\infty}\myint{\prt\Gw_{n}}{}
\abs{Z_{f,n}-Z_{0,n}}dS\,dt=0
$$
by approximations. By the maximum principle $Z_{f,n,k}-Z_{0,n,k}\leq 0$. Letting $n\to\infty$ yields to
$$\overline u_{f}-u_{f}\leq \overline u_{0}-u_{0},
$$
which ends the proof.\qeda
\section{The case $1<q<N/(N-2)$}
\setcounter{equation}{0}
In this section we assume that $\Gw$ is a domain of $\BBR^N$ with a compact boundary. We first prove that the maximal solution is a large solution

\bth{max=large} Assume $1<q<N/(N-2)$ and $f\in L^1_{loc\,+}(\Gw)$
. Then the maximal solution $\overline u_{f}$ of (\ref{E1}) in $Q_{T}^\Gw$ which satisfies (\ref{C1}) satisfies also (\ref{C4}).
\es
\Proof In Appendix we construct the self-similar solution $V:=V_{N}$ of
(\ref{E1}) in $Q_{\infty}^{\BBR^N\setminus\{0\}}$ which has initial trace zero in $\BBR^N\setminus\{0\}$ and satisfies
$$\lim_{|x|\to 0}V_{N}(x,t)=\infty,$$ 
locally uniformly on $[\gt,\infty)$, for any $\gt>0$. Furthermore
$V_{N}(x,t)=t^{-1/(q-1)}H_{N}(|x|/\sqrt t)$.
If $a\in\prt\Gw$, the restriction to $\Gw_{n}$ of the function $V_{N}(x-a,t)$ is bounded from above by $u_{n,f}$. Letting $n\to\infty$ yields to 
\begin{equation}\label{V5}
V_{N}(x-a,t)\leq \overline u_{f}(x,t)\quad\forall (x,t)\in Q_{\infty}^\Gw.
\end{equation}
If we consider $x\in\Gw$ and denote by $a_{x}$ a projection of $x$ onto $\prt\Gw$, there holds 
\begin{equation}\label{V6}
t^{-1/(q-1)}H_{N}(\gr(x)/\sqrt t)=V_{N}(x-a_{x},t)\leq \overline u_{f}(x,t).
\end{equation}
Using (\ref{V2}), we derive that $\overline u_f$ satisfies (\ref{C4}).\qeda

\bth{UniqueI} Assume $1<q<N/(N-2)$, $f\in L^1_{loc\,+}(\Gw)$ and $\prt\Gw=\prt\overline\Gw^c$. Then $\overline u_{f}$ is the unique solution of (\ref{E1}) in $Q_{T}^\Gw$ which satisfies (\ref{C1}) and (\ref{C4}).
\es
\Proof Assume that $u_{f}$ is a solution of (\ref{E1}) in $Q_{T}^\Gw$ such that (\ref{C1}) and (\ref{C4}) hold. By \rth{contraction} there exists a positive solution $u_{0}$ with zero initial trace such that 
\begin{equation}\label{UK}
0\leq u_{f}- u_{0}\leq\underline u_{f}
\end{equation}
and (\ref{D2}) are satisfied. Since $\underline u_{f}(x,t)\leq \left((q-1)t\right)^{-1/(q-1)}$ (notice that this last expression is the maximal solution of (\ref{E1}) in $Q_{\infty}^{\BBR^N}$), the function $u_{0}$ satisfies also (\ref{C4}). Therefore, it is sufficient to prove 
that $\overline u_{0}=u_{0}:=u$.\smallskip

\noindent{\it Step 1: bilateral estimates.} Since $\prt\Gw=\prt\overline\Gw^c$, for any $a\in\prt\Gw$, there exists a sequence $\{a_{n}\}\subset\overline\Gw^c$ converging to $a$. If $u$ is any solution of (\ref{E1}) in $Q_{T}^\Gw$ which satisfies (\ref{C4}) and (\ref{D0}), there holds
$$V_{N}(x-a_{n},t)\leq u(x,t)\Longrightarrow V_{N}(x-a,t)\leq u(x,t).
$$
In particular, if $a=a_{x}$, we see that $u$ satisfies (\ref{V6}). In order to obtain an estimate from above we consider for $r<\gr(x)$ the solution 
$(y,t)\mapsto u_{x,r}(y,t)$ of 
\begin{equation}\label{V7}\left\{\BA {l}
\prt_{t}u_{x,r}-\Gd u_{x,r}+u_{x,r}^q=0\quad\text{in }\;Q_{\infty}^{B_{r}(x)}\\[2mm]
\lim_{(y,t)\to (z,0)}u_{x,r}(y,t)=0\quad\forall z\in B_{r}(x)\\[2mm]
\lim_{|x|\uparrow r}u_{x,r}(x,t)=\infty\quad\text{locally uniformly on $[\gt,\infty)$, 
for any $\gt>0$}\EA\right. 
\end{equation}
Then 
$$\overline u_{0} (y,t)\leq u_{x,r}(y,t)\Longrightarrow \overline u_{0} (y,t)\leq u_{x,\gr(x)}(y,t)\quad\forall (y,t)\in Q_{\infty}^{B_{\gr(x)}(x)}. $$
In particular, with $u_{0,r}=u_{r}$,
$$\overline u_{0} (x,t)\leq u_{\gr(x)}(0,t)=(\gr(x))^{-2/(q-1)}
u_{1}(0,t/(\gr(x))^{2}).$$
Therefore
\begin{equation}\label{V8}
t^{-1/(q-1)}H_{N}(\gr(x)/\sqrt t)\leq u(x,t)\leq \overline u_{0}(x,t)\leq (\gr(x))^{-2/(q-1)}
u_{1}(0,t/(\gr(x))^{2}).
\end{equation}
The function $s\mapsto u_{1}(0,s)$ is increasing by the same argument as the one of \rcor{sigma=0} and bounded from above by the unique solution $P$ of 
\begin{equation}\label{V9}\left\{\BA {l}
-\Gd P+P^q=0\quad\text{in } B_{1}\\
\lim_{|x|\to 1}P(x)=\infty.
\EA\right. 
\end{equation}
Therefore it converges to $P$ locally uniformly in $B_{1}$ and $\displaystyle\lim_{s\to\infty}u_{1}(0,s)= P(0)$. Thus 
\begin{equation}\label{V9'}t/(\gr(x))^{2}\to\infty\Longrightarrow (\gr(x))^{-2/(q-1)}
u_{1}(0,t/(\gr(x))^{2})\approx P(0)(\gr(x))^{-2/(q-1)}.
\end{equation}
On the other hand, if $t/(\gr(x))^{2}\to\infty$, equivalently $\gr(x)/\sqrt t\to 0$,
\begin{equation}\label{V9''} t^{-1/(q-1)}H_{N}(\gr(x)/\sqrt t)\approx \gl_{N,q}t^{-1/(q-1)}(\gr(x)/\sqrt t)^{-2/(q-1)}=\gl_{N,q}(\gr(x))^{-2/(q-1)},
\end{equation}
by (\ref{V4}).\smallskip

Next, in order to obtain an estimate from above of $u_{1}(0,s)$ when $s\to 0$, we compare $u_{1}$ to a solution $u_{\Gth}$ of (\ref{E1}) in $Q_{\infty}^\Gth$, where $\Gth$ is a polyhedra inscribed in $B_{1}$; this polyhedra is a finite intersection of half spaces $\Gg_{i}$  containing $\Gp$. In each of the half space $\Gg_{i}$, with boundary $\gg_{i}$, we can consider the solution $W_{i}$ of (\ref{E1}) in $Q_{\infty}^{\Gg_{i}}$ which tends to infinity on $\gg_{i}\ti (0,\infty)$ and has value $0$ on 
$\Gg_{i}\ti\{0\}$. This solution depends only on the distance to $\gg_{i}$ and $t$. Thus it is expressed by the function $V_{1}$ defined in \rprop{singsol1} when $N=1$. Moreover, since a sum of solutions is a super solution, 
\begin{equation}\label{V10}u_{1}\leq  u_{\Gth}\leq\sum_{i}W_{i}\Longrightarrow
u_{1}(0,s)\leq \sum_{i}H_{1}(\dist (0,\gg_{i})/\sqrt s).
\end{equation}
We can choose the hyperplanes $\gg_{i}$ such that for any $\gd\in (0,1)$, there exists $C_{\gd}\in \BBN_{*}$ such that 
\begin{equation}\label{V11}u_{1}(0,s)\leq C_{\gd}H_{1}( (1-\gd)/\sqrt s).
\end{equation}
Using (\ref{V3}) we derive
$$
u(x,t)\geq c_{N,q}(\gr(x))^{2/(q-1)-N}t^{N/2-1/(q-1)}e^{-(\gr(x))^2/4 t},
$$
when $\gr(x)/\sqrt t\to\infty$, and 
$$\overline u_{0}(x,t)\leq CH_{1}((1-\gd)\gr(x)/\sqrt t)
\leq C(1-\gd)^{2/(q-1)-1}(\gr(x))^{2/(q-1)-1}t^{1/2-1/(q-1)}e^{-((1-\gd)\gr(x))^2/4 t}.
$$ 
Therefore, there exists $\gth>1$ such that
\begin{equation}\label{V11'}\overline u_{0}(x,t)\leq C(\gr(x))^{2/(q-1)-N}t^{N/2-1/(q-1)}e^{-(\gr(x))^2/4 \gth t}\leq Cu(x,\gth t),
\end{equation}
when $\gr(x)/\sqrt t\to\infty$. Finally, when $m^{-1}\leq\gr(x)/\sqrt t\leq m$ for some $m>1$, (\ref{V8}) shows that $ (\gr(x))^{-2/(q-1)}
u_{1}(0,t/(\gr(x))^{2})$ and $t^{-1/(q-1)}H_{N}(\gr(x)/\sqrt t)$ are comparable. In conclusion, there exist constants $C>P(0)/\gl_{N,q}>1$ and $\gth>1$ such that 
\begin{equation}\label{V12}
u(x,t)\leq \overline u_{0}(x,t)\leq Cu(x,\gth t)\quad\forall (x,t)\in Q_{\infty}^{\Gw}.
\end{equation}

\noindent{\it Step 2: End of the proof.} Let $\gt>0$ and $C'>C$ be fixed. The function 
$$t\mapsto u_{\gt}(x,t):=C'u(x,t+\gth\gt)$$ 
is a supersolution of (\ref{E1}) in $\Gw\ti (0,\infty)$ which satisfies $u_{\gt}(x,0)=C'u(x,\gth\gt)> \overline u_{0}(x,\gt)$ by  (\ref{V12}).
Furthermore,
$$C'u(x,t+\gth\gt)\geq C'(t+\gth\gt)^{_{-1/(q-1)}}H_{N}(\gr(x)/\sqrt{t+\gth\gt})=C'\gl_{N,q}(1+o(1))(\gr(x))^{-2/(q-1)},
$$
as $\gr(x)\to 0$, locally uniformly for $t\in [0,\infty)$. Similarly, 
$$\overline u_{0}(x,t+\gt)\leq (\gr(x))^{-2/(q-1)}
u_{1}(0,(t+\gt)/(\gr(x))^{2})=P(0)(1+o(1))(\gr(x))^{-2/(q-1)},
$$
as $\gr(x)\to 0$, and also locally uniformly for $t\in [0,\infty)$. Therefore
 $(\overline u_{0}(x,t)-u_{\gt}(x,t))_{+}$ vanishes in a neighborhood of $\prt\Gw\ti [0,T]$ for any $T>0$. By the maximum principle
$$u_{\gt}(x,t)\geq \overline u_{0}(x,t)\quad\forall (x,t)\in \Gw\ti (0,\infty).
$$
Letting $\gt\to 0$ and $C'\to C$ yields to
\begin{equation}\label{V13}
u(x,t)\leq \overline u_{0}(x,t)\leq Cu(x, t)\quad\forall (x,t)\in Q_{\infty}^{\Gw}.
\end{equation}
The conclusion of the proof is contradiction, following an idea introduced in \cite{MV2} and developped by \cite{Ve1} in the elliptic case. We assume $u\neq\overline u_{0}$, thus $u<\overline u_{0}$. By convexity the function
$$w=u-\myfrac{1}{2C}(\overline u_{0}-u)
$$
is a supersolution and $w<u$. Moreover $w>w':=((1+C)/2C)u$ and $w'$ is a subsolution. Consequently, there exists a solution $u_{1}$ of (\ref{E1}) which satisfies
\begin{equation}\label{V14}
w'<u_{1}\leq w\Longrightarrow \overline u_{0}-u_{1}\geq \left(1+K^{-1}\right) (\overline u_{0}-u)\quad\text{ in }\, Q_{\infty}^{\Gw}.
\end{equation}
Notice that $u_{1}$ satisfies (\ref{D0}) and (\ref{C4}), therefore it satisfies (\ref{V13}) as $u$ does it. Replacing $u$ by $u_{1}$ and introducing the supersolution 
$$w_{1}=u_{1}-\myfrac{1}{2C}(\overline u_{0}-u_{1})
$$
and the subsolution $w_{1}':=((1+C)/2C)u_{1}$ we see that there exists a solution $u_{2}$  of (\ref{E1}) such that
\begin{equation}\label{V15}
w_{1}'<u_{2}\leq w_{1}\Longrightarrow \overline u_{0}-u_{2}\geq \left(1+K^{-1}\right)^2 (\overline u_{0}-u)\quad\text{ in }\, Q_{\infty}^{\Gw}.
\end{equation}
By induction, we construct a sequence of positive solutions $u_{k}$ of (\ref{E1}), subject to (\ref{D0}) and (\ref{C4}) such that 
\begin{equation}\label{V16}
\overline u_{0}-u_{k}\geq \left(1+K^{-1}\right)^k (\overline u_{0}-u)\quad\text{ in }\, Q_{\infty}^{\Gw}.
\end{equation}
This is clearly a contradiction since $\left(1+K^{-1}\right)^k\to\infty$ as $k\to\infty$ and $\overline u_{0}$ is locally bounded in $Q_{\infty}^{\Gw}$.
\qeda
\section{The local continuous graph property}
\setcounter{equation}{0}
In this section, we assume that $\prt\Gw$ is compact and is locally the graph of a continuous function, which means that there exists a finite number of open sets $\Gw_{j}$ ($j=1,...,k$) such that $\Gg\cap \Gw_{j}$ is the graph  of a continuous function. Our main result is the following
\bth{uniqth} Assume $q>1$ and $f\in L^1_{loc\,+}(\Gw)$.  Then there exists at most one positive solution of (\ref{E1}) in $Q_{\infty}^\Gw$ satisfying (\ref{C1}) and (\ref{C4}).
\es

Suppose $u_{f}$ satisfies (\ref{E1}) in $Q_{\infty}^\Gw$ satisfying (\ref{C1}) and (\ref{C4}), then clearly the maximal solution $\overline u_{f}$ endows the same properties. In order to prove that $u_{f}=\overline u_{f}$, we can assume that $f=0$ by \rth{contraction}. We denote by $u$ this large solution with zero initial trace. We consider some $j\in \{1,...,k\}$, perform a rotation,  denote by $x=(x',x_{N})\in \BBR^{N-1}\ti\BBR$ the coordinates in $\BBR^N$ and represent $\Gg\cap \Gw_{j}$ as the graph of a continuous positive function $\gf$ defined in $C=\{x'\in\BBR^{N-1}:|x'|\leq R\}$. We identify $C$ with $\{x=(x',0):|x'|\leq R\}$ and set
$$\Gg_{1}=\{x=(x',\gf(x')):x'\in C\},
$$
$$\Gg_{2}=\{x=(x',x_{N}):x'\in \prt C,\,0\leq x_{N}< \gf(x'), \},
$$
and 
$$G_{R}=\{x\in\BBR^N:|x'|<R,\,0<x_{N}<\gf(x')\}.
$$
We can assume that $\overline G_{R}\subset \Gw\cup\Gg_{1}$,
$$\inf\{\gf(x'):x'\in C\}=R_{0}>0\quad\text{and }\;
\sup\{\gf(x'):x'\in C\}=R_{1}>R_{0}.$$
For $\gs> 0$, small  enough, we consider 
$\gf_{\gs}\in C^\infty(C)$ satisfying
$$\gf(x')-\gs/2\leq \gf_{\gs}(x')\leq \gf(x')+\gs/2\quad\forall x'\in C,
$$
and set 
$$G_{\gs,R}=\{x\in\BBR^N:|x'|<R,\,0<x_{N}<\gf_{\gs}(x')-\gs\}
$$
and
$$G'_{\gs,R}=\{x\in\BBR^N:|x'|<R,\,0<x_{N}<\gf_{\gs}(x')+\gs\}.
$$
The upper boundaries of $G_{\gs}$ and $G'_{\gs}$ are defined by
$$\Gg_{1,\gs}=\{x=(x',\gf_{\gs}(x')-\gs):x'\in C\},
$$
$$\Gg'_{1,\gs}=\{x=(x',\gf_{\gs}(x')+\gs):x'\in C\},
$$
and the remaining boundaries are
$$\Gg_{2,\gs}=\{x=(x',x_{N}):x'\in \prt C,\,0\leq x_{N}\leq \gf_{\gs}(x')-\gs \},
$$
$$\Gg'_{2,\gs}=\{x=(x',x_{N}):x'\in \prt C,\,0\leq x_{N}\leq \gf_{\gs}(x')+\gs \}.
$$
In order to have the monotonicity of the domains, we can also assume
\begin{equation}\label{L1}
\gf_{\gs}(x')-\gs<\gf_{\gs'}(x')-\gs'<\gf_{\gs'}(x')+\gs'<
\gf_{\gs}(x')+\gs\quad\forall \,0<\gs'<\gs\quad\forall \,x'\in C,
\end{equation}
thus, under the condition $0<\gs'<\gs$,
\begin{equation}\label{L2}
G_{\gs,R}\subset G_{\gs',R}\subset G_{R}\subset G'_{\gs',R}\subset G'_{\gs,R}.
\end{equation}
The localization procedure is to consider the restriction of $u$ to $Q_{\infty}^{G_{R}}:=G_{R}\ti (0,\infty)$, thus $u$ is regular in $G_{R}\cup\Gg_{2}\ti [0,\infty)$ and satifies
\begin{equation}\label{C5}
\lim_{x_{N}\to\gf(x')}u(x',x_{N},t)=\infty,
\end{equation}
uniformly with respect to $(x',t)\in C\ti[\gt,T]$, for any $0<\gt<T$. We construct $v_{\gs}$ as solution of
\begin{equation}\label{E2}
\prt_{t}v_{\gs}-\Gd v_{\gs}+v_{\gs}^q=0\quad\text{in }
Q_{\infty}^{G_{\gs,R}}:=G_{\gs,R}\ti (0,\infty),
\end{equation}
subject to the initial condition
\begin{equation}\label{v1}
\lim_{t\to 0}v_{\gs}(x,t)=0\quad\text{locally uniformly in } G_{\gs,R},
\end{equation}
and the boundary conditions
\begin{equation}\label{v2}
\lim_{x_{N}\to \gf_{\gs}(x')-\gs}v_{\gs}(x',x_{N},t)=\infty\quad\forall (x',t)\in C\ti (0,\infty],
\end{equation}
uniformly on any set $K\ti [\gt,T]$, where $T>\gt>0$ and $K$ is a compact subset of $C$ and
\begin{equation}\label{v3}
v_{\gs}(x,t)=0\quad\forall (x,t)\in \Gg_{2,\gs}\ti [0,\infty).
\end{equation}
We also construct $w_{\gs}$ as solution of
\begin{equation}\label{E'2}
\prt_{t}w_{\gs}-\Gd w_{\gs}+w_{\gs}^q=0\quad\text{in }
Q_{T}^{G'_{\gs,R}}:=G'_{\gs,R}\ti (0,\infty),
\end{equation}
subject to the initial condition
\begin{equation}\label{w1}
\lim_{t\to 0}w_{\gs}(x,t)=0\quad\text{locally uniformly in } G'_{\gs,R}, 
\end{equation}
and the boundary conditions
\begin{equation}\label{w2}\left\{\BA {l}
(i)\qquad w_{\gs}(x,t)=0\quad\forall (x,t)\in \Gg'_{1,\gs}\ti [0,T],\\[2mm]
(i')\qquad \lim_{(x,s)\to (y,t)}w_{\gs}(x,t)=\infty\quad\forall (y,s)\in \Gg'_{2,\gs}\ti [0,T].
\EA\right.\end{equation}

The functions $v_{\gs}$ and $w_{\gs}$ inherit the following properties in which the local graph property plays a fundamental role, allowing translations of the truncated domains in the $x_{N}$-direction.

\blemma{prop}For $\gs>\gs'>0$ there holds
\begin{equation}\label{P1}
v_{\gs'}\leq v_{\gs}\quad \text{in }Q_{\infty}^{G_{\gs,R}},
\end{equation}
\begin{equation}\label{P2}
w_{\gs'}\leq w_{\gs}\quad \text{in }Q_{\infty}^{G'_{\gs',R}},
\end{equation}
\begin{equation}\label{P3}\left.\BA{l}
(i)\quad v_{\gs}(x',x_{N}-2\gs,t)\leq u(x',x_{N},t)\quad \text{in }Q_{\infty}^{ G_{R}}\\[2mm]
(ii)\quad u(x',x_{N},t)\leq v_{\gs}(x,t)+w_{\gs}(x,t)\quad \text{in }Q_{\infty}^{G_{\gs,R}}.
\EA\right.\end{equation}
\es
\Proof The inequalities (\ref{P1}) and (\ref{P2}) are the direct consequence of the fact that the domains $G_{\gs,R}$ and $G'_{\gs',R}$ are Lipschitz and the functions $v_{\gs}$ and $w_{\gs}$ are constructed by approximations of solutions of (\ref{E1}) with bounded boundary data. For proving (\ref{P3})-(i), we compare, for $\gt>0$, 
$u(x,t-\gt)$ and $v_{\gs}(x',x_{N}-2\gs,t)$ in $Q_{\infty}^{G_{R}}$. Because $u$ satisfies (\ref{C4}), and $v_{\gs}(x',x_{N}-2\gs,0)=0$ in $G_{R}$, (\ref{P3})-(i) follows by the maximum principle. The proof of (\ref{P3})-(ii) needs no translation, but the fact that the sum of two solutions is a supersolution.\qeda
\bcor{sigma=0}There exist $v_{0}=\displaystyle\lim_{\gs\to 0}v_{\gs}$ and $w_{0}=\displaystyle\lim_{\gs\to 0}w_{\gs}$ and there holds
\begin{equation}\label{P4}
v_{0}\leq u\leq v_{0}+w_{0}\quad \text{in }Q_{\infty}^{ G_{R}}.
\end{equation}
Moreover, the functions $t\mapsto v_{0}(x,t)$ and  $t\mapsto w_{0}(x,t)$ are increasing on $(0,\infty)$, $\forall x\in G_{R}$.
\es
\Proof The first assertion follows from (\ref{P1})-(\ref{P2}), and (\ref{P4}) from (\ref{P3}). Since $v_{0}$ is the limit, when $\gs\to 0$ of $v_{\gs}$ which satisfy equation (\ref{E2}) in $Q_{T}^{G_{\gs,R}}$, initial condition (\ref{v1}) and boundary conditions (\ref{v2}), (\ref{v3}), it is sufficient to prove the monotonicity of  $t\mapsto v_{\gs}(.,t)$. Moreover  $v_{\gs}$ is the limit, when $k$ tends to infinity of the $v_{k,\gs}$ solutions of (\ref{E1}) in $Q_{T}^{G_{\gs,R}}$, which satisfy the same boundary conditions as $v_{\gs}$ on $\Gg_{2,\gs}\ti [0,T]$, the same zero initial condition and 
$$\lim_{x_{N}\to\gf(x')-\gs}v_{k,\gs}(x',x_{N},t)=k.
$$
For $\gt>0$, we define $V_{\gt}$  by  $V_{\gt}(x,t)=(v_{k,\gs}(x,t)-v_{k,\gs}(x,t+\gt))_{+}$. Because $\prt G_{\gs,R}$ is Lipschitz and $V_{\gt}$ is a subsolution of (\ref{E1})  which vanishes on $\prt G_{\gs,R}\ti [0,T]$ and at $t=0$, it is identically zero. This implies $v_{k,\gs}(x,t)\leq v_{k,\gs}(x,t+\gt)$, and the monotonicity property of $v_{0}$, by strict maximum principle and letting $\gs\to 0$. The proof of the 
monotonicity of $w_{0}$ is similar.\qeda\medskip

The key step of the proof is the following result.
\bprop{lim} Let $\ge,\gt>0$. Then there exists $\gd_{\ge}>0$ such that, if we denote
$$ G_{\gd,R'}=\{x=(x',x_{N}):|x'|<R'\text{ and }\gf(x')-\gd \leq x_{N}<\gf(x')\},
$$
there holds, for $R'<R/\sqrt{N-1}$,
\begin{equation}\label{Q1}w_{0}(x,t)\leq \ge v_{0}(x,t+\gt)\quad \forall (x,t)\in Q_{\infty}^{G_{\gd,R'}}.
\end{equation}
\es
\Proof Using the result in Appendix, we recall that $V:=V_{1}$ is the  unique positive and self-similar solution of the problem
\begin{equation}\label{W1}\left\{\BA {l}
\prt_{t}V-\prt_{zz}V+V^q=0\quad\text{in }\;\BBR_{+}\ti\BBR_{+}\\[2mm]
\lim_{t\to 0}V(z,t)=0\quad\forall z>0\\[2mm]
\lim_{z\to 0}V(z,t)=\infty\quad\forall t>0,
\EA\right. 
\end{equation}
and it is expressed by $V_{1}(z,t)=t^{-1/(q-1)}H_{1}(x/\sqrt t)$, where $H_{1}$ satisfies (\ref{V2})-(\ref{V3}) with $N=1$. We set 
$R_{N}=R/\sqrt {N-1}$ so that 
$$C_{\infty}:=\{x'=(x_{1},...,x_{N-1}):\sup_{j\leq N-1}{|x_{j}|}<R_{N}\}\subset C=\{x':|x'|\leq R\}$$ 
and we define
$$\tilde w(x,t)=W(x_{N},t)+\sum_{j=1}^{N-1}(W(x_{j}-R,t)+W(R-x_{j},t)).
$$
The function $\tilde w$ a super solution in $\Gth\ti \BBR^+$ where $\Gth:=\{(x',x_{N}):x'\in C_{\infty},x_{N}>0\}$ which blows up on 
$$\{x:x_{N}=0\,,\;\sup_{j\leq N-1}{|x_{j}|}\leq R \}\bigcup_{j\leq N-1}
\left\{x:x_{N}\geq 0\,,\;x_{j}=\pm R\right\}.
$$
Therefore $w_{0}\leq \tilde w$ in $Q_{T}^{G_{R_{N}}}$. Moreover
$\tilde w(x,t)\to 0$ when $t\to 0$, uniformly on 
$$G^*_{\ga,R'}:=\{x=(x_{1},x_{2}):|x_{1}|\leq R',\ga\leq x_{2}\leq\gf(x_{1})\},$$
for any $\ga\in (0,R_{0}]$ and $R'\in (0,R_{N})$. Since for any $\gt>0$, $v_{0}(x,t+\gt)\to \infty$ when $\gr(x)\to 0$, locally uniformly on $[0,\infty)$, and $\tilde w(x,t)$ remains uniformly bounded on $Q_{\infty}^{G_{\gd,R'}}$, for any $\gd>R_{0}$, it follows that for any $\ge>0$ there exists $\gd_{\ge}>0$ such that 
$$w_{0}(x,t)\leq \tilde w(x,t)\leq \ge v_{0}(x,t+\gt)\quad\forall (x,t)\in 
Q_{\infty}^{G_{\gd_{\ge},R'}}.$$\qeda\medskip

\noindent{\it Proof of \rth{uniqth}.} Assume $u$ is a solution of (\ref{E1}) satisfying (\ref{C1}) and (\ref{C4}). Then there holds in
$Q_{\infty}^{G_{\gd_{\ge},R'}}$,
\begin{equation}\label{W5}
v_{0}(.,t)\leq u (.,t)\leq v_{0}(.,t)+\ge v_{0}(.,t+\gt).
\end{equation}
Therefore
$$v_{0}(.,t+\gt)\leq u (.,t+\gt)\leq v_{0}(.,t+\gt)+\ge v_{0}(.,t+2\gt),
$$
from which follows
$$ (1+\ge)u (.,t+\gt)\geq (1+\ge)v_{0}(.,t+\gt)\geq v_{0}(.,t)+
\ge v_{0}(.,t+\gt)
$$
since $t\mapsto v_{0}(.,t)$ is increasing by \rcor{sigma=0}. The maximal solution $\overline u_{0}$ satisfies (\ref{W5}) too; consequently the following inequality is verified in $Q_{\infty}^{G_{\gd_{\ge},R'}}$,
\begin{equation}\label{W6}(1+\ge)u (.,t+\gt)\geq \overline u_{0}(.,t).
\end{equation}
Since $\prt\Gw$ is compact, there exists $\gd^*>0$ such that (\ref{W6}) holds
whenever $t\in [0,T]$ ($T>0$ arbitrary) and $\gr(x)\leq \gd^*$. Furthermore 
$$\lim_{t\to 0} \max\{(\overline u_{0}(x,t)-(1+\ge)u (x,t+\gt))_{+}:\gr(x)\geq\gd^*\}=0
$$
because of (\ref{C1}). Since $(\overline u_{0}(x,t)-(1+\ge)u (x,t+\gt))_{+}$ is a subsolution, which vanishes at $t=0$ and near $\prt\Gw\ti [0,T]$, it follows that (\ref{W6}) holds in $Q_{T}^\Gw$. Letting $\ge\to 0$ and $\gt\to 0$ yields to $u\geq \overline u_{0}$. \qeda\medskip

\noindent\Remark The existence of large solutions when $q\geq N/(N-2)$ is a difficult problem as it is already in the elliptic case. We conjecture  that the necessary and sufficient conditions, obtained by Dhersin-Le Gall when $q=2$ \cite {DL} and Labutin \cite {La} in the general case $q>1$, and expressed by mean of a Wiener type criterion involving the $C_{2,q'}^{\BBR^N}$-Bessel capacity, are still valid. As in \cite{MV1}, it is clear that if $\prt\Gw$ satisfies the exterior segment property and $1<q<(N-1)/(N-3)$, then $\overline u_{0}$ is a large solution.

\section{Appendix}
\setcounter{equation}{0}
The proof of this  result is based upon the existence of solution of (\ref{E1}) in $Q_{\infty}^{\BBR^N\setminus\{0\}}$ with a persistent singularity on $\{0\}\ti [0,\infty)$.
\bprop{singsol1}For any $q>1$, there exists a unique positive function $V:=V_{N}$ defined in $\BBR_{+}\ti\BBR_{+}$ satisfying, for any $\gt>0$
\begin{equation}\label{V1}\left\{\BA {l}
\prt_{t}V-\Gd V+V^q=0\quad\text{in }\;Q_{\infty}^{\BBR^N\setminus\{0\}}\\[2mm]
\lim_{(x,t)\to (y,0)}V(x,t)=0\quad\forall y\in \BBR^N\setminus\{0\}\\[2mm]
\lim_{|x|\to 0}V(x,t)=\infty\quad\text{locally uniformly on $[\gt,\infty)$, 
for any $\gt>0$}\EA\right. 
\end{equation}
Then $V_{N}(x,t)=t^{-1/(q-1)}H_{N}(|x|/\sqrt t)$, where $H:=H_{N}$ is the unique positive function satisfying
\begin{equation}\label{V2}\left\{\BA {l}
H''+\left(\myfrac{N-1}{r}+\myfrac{r}{2}\right)H'+\myfrac{1}{q-1}H-H^q=0\quad\text{in }\;\BBR_{+}\\[3mm]
\lim_{r\to 0}H(r)=\infty\\[2mm]
\lim_{r\to \infty}r^{2/(q-1)}H(r)=0.
\EA\right. 
\end{equation}
Furthermore there holds
\begin{equation}\label{V3}
H_{N}(r)=c_{N,q}r^{2/(q-1)-N}e^{-r^2/4}(1+O(r^{-2}))\quad\text{as }\;r\to\infty,
\end{equation}
and
\begin{equation}\label{V4}
H_{N}(r)=\gl_{N,q}r^{-2/(q-1)}(1+O(r))\quad\text{as }\;r\to 0,
\end{equation}
\es
\Proof If we assume $1<q<N/(N-2)$, the $C_{2,1,q'}$ parabolic capacity of the axis $\{0\}\ti\BBR\subset \BBR^{N+1}$ is positive, therefore there exists a unique solution $u:=u_{\gm}$ to the problem
\begin{equation}\label{A1}\BA {l}
\prt_{t}u-\Gd u+|u|^{q-1}u=\gm\quad\in \BBR^N\ti\BBR,
\EA\end{equation}
(see \cite{BP}) where $\gm$ is the uniform measure on $\{0\}\ti\BBR_{+}$ defined by
$$\myint{}{}\gz d\gm=\myint{0}{\infty}\gz(0,t)dt\quad\forall\gz\in C^\infty_{0}(\BBR^{N+1}).
$$
If we denote $T_{\ell}[u](x,t)=\ell^{2/(q-1)} u(\ell x,\ell^2 t)$ for $\ell>0$, then $T_{\ell}$ leaves the equation (\ref{E1}) invariant, and
$T_{\ell}[u_{\gm}]=u_{\ell^{2/(q-1)-N}\gm}$. If we replace $\gm$ by $k\gm$ ($k>0$), we obtain
\begin{equation}\label{A2}
T_{\ell}[u_{k\gm}]=u_{\ell^{2/(q-1)-N}k\gm}.
\end{equation}
Moreover, any solution of (\ref{E1}) in $\BBR^N\setminus\{0\}\ti\BBR_{+}$ which vanishes on $\BBR^N\setminus\{0\}\ti\{0\}$ is bounded from above by the maximum solution $u:=U$ of 
\begin{equation}\label{A3}
-\Gd u+u^q=0\quad\text{in }\BBR^N\setminus\{0\}.
\end{equation}
This is obtained by considering the solution $U_{\ge}$ of 
\begin{equation}\label{A4}\left\{\BA {l}
-\Gd u+u^q=0\quad\text{in }\BBR^N\setminus \overline B_{\ge}\\[2mm]
\,\displaystyle\lim_{|x|\to\ge}u(x)=\infty.
\EA\right.\end{equation}
Actually, 
\begin{equation}\label{A5}
U(x):=\lim_{\ge\to 0}U_{\ge}(x)=\gl_{N,q}|x|^{-2/(q-1)}\quad\text{with }
\gl_{N,q}:=\left[\left(\myfrac{2}{q-1}\right)\left(\myfrac{2q}{q-1}-N\right)\right]^{1/(q-1)},
\end{equation}
an expression which exists since $1<q<N/(N-2)$. If we let $k\to\infty$ in (\ref{A2}), using the monotonicity of $\gm\mapsto u_{\gm}$, we obtain that $u_{k\gm}\to u_{\infty\gm}$, $u_{\infty\gm}\leq U$ and
\begin{equation}\label{A6}
T_{\ell}[u_{\infty\gm}]=u_{\ell^{2/(q-1)-N}\infty\gm}=u_{\infty\gm}\quad\forall\ell>0.
\end{equation}
This implies that $u_{\infty\gm}$ is self-similar, that is
$$u_{\infty\gm}(x,t)=t^{-1/(q-1)}h(x/\sqrt t).
$$
Furthermore, $h(.)$ is positive and radial as $x\mapsto u_{\gm}(x,t)$ is, and it solves
\begin{equation}\label{A7}
h''+\left(\myfrac{N-1}{r}+\myfrac{r}{2}\right)h'+\myfrac{1}{q-1}h-h^q=0\quad\text{in }\;\BBR_{+}.
\end{equation}
Since $u_{\gm}(x,0)=0$ for $x\neq 0$, the a priori bounds $u_{k\gm}\leq U$, the equicontinuity of the $\{u_{k\gm}\}_{k>0}$ implies that
$u_{\infty\gm}(x,0)=0$ for $x\neq 0$; therefore
\begin{equation}\label{A8}
\lim_{r\to\infty}r^{2/(q-1)}h(r)=0.
\end{equation}
The same argument as the one used in the proof of \rcor{sigma=0} implies that $t\mapsto u_{\gm}(x,t)$ is increasing, therefore $\lim_{x\to 0}u_{\gm}(x,t)=\infty$ for $t>0$. This implies 
 $\lim_{r\to 0}h(r)=\infty$. Then the proof of (\ref{V3}) follows from \cite[Appendix]{MV4}. When $r\to 0$, $h$ could have two possible behaviours \cite{Ve2}:\smallskip
 
 \noindent (i) either
 \begin{equation}\label{A9}
h(r)=\gl_{N,q}r^{-2/(q-1)}(1+O(r)),
\end{equation}

 \noindent (ii) or there exists $c\geq 0$ such that
  \begin{equation}\label{A10}
h(r)=cm_{N}(r)(1+O(r)),
\end{equation}
where $m_{N}(r)$ is the Newtonian kernel if $N\geq 2$ and $m_{1}(r)=1+o(1)$.\smallskip

If (ii) were true with $c>0$ (the case $c=0$ implying that $h=0$ because of the behavior at $\infty$ and maximum principle), it would lead to
  \begin{equation}\label{A10'}
u_{\infty\gm}(x)=c|x|^{2-N}t^{N-2-1/(q-1)}(1+o(1))\quad\text{as }x\to 0,
\end{equation}
for all $t>0$. Therefore
  \begin{equation}\label{A11}\myint{\ge}{T}\myint{B_{1}}{}u^q_{k\gm}dx\,dt<C(\ge),
\end{equation}
for any $\ge>0$ and $k\in (0,\infty]$. We write (\ref{A1}) under the form
$$\prt_{t}u_{k\gm}-\Gd u_{k\gm}=g_k+k\gm
$$
where $g_{k}=-u^q_{k\gm}$, then $u_{k\gm}=u'_{k\gm}+u_{k}''$, where 
$$\prt_{t}u'_{k\gm}-\Gd u'_{k\gm}=k\gm
$$
and 
$$\prt_{t}u''_{k}-\Gd u''_{k}=g_{k}.
$$
By linearity $u'_{k\gm}=ku'_{\gm}$. Because of (\ref{A11}) $u''_{k}$ remains uniformly bounded in $L^{1}(B_{1}\ti (\ge,T)$. This clearly contradicts $\lim_{k\to\infty}u'_{k\gm}=\infty$. Thus (\ref{V4}) holds.
The proof of uniqueness is an easy adaptation of \cite[Lemma 1.1]{MV1}: the fact that the domain is not bounded being compensated by the strong decay estimate (\ref{V3}). This unique solution is denoted by $V_{N}$ and $h=H_{N}$.\qeda
 
 \end{document}